\documentclass[12pt,a4paper]{article}
\usepackage{mathtext}
\usepackage[cp1251]{inputenc}
\usepackage[english]{babel}
\usepackage{euscript,amssymb,amsfonts,amsmath,amsthm}
\usepackage{mathrsfs}
\usepackage[T2A]{fontenc}
\usepackage{color}

\newtheorem{theorem}{\indent \sc Theorem}

\newtheorem{cor}{\indent \sc Corollary}

\theoremstyle{remark}

\newcommand{\eps}{\varepsilon}
\newcommand{\R}{\mathbb R}
\newcommand{\N}{\mathbb N}
\newcommand{\E}{{\sf E}}

\newcommand{\Prob}{{\sf P}}

\renewcommand{\le}{\leqslant}
\renewcommand{\ge}{\geqslant}
\newcommand{\eqd}{\stackrel{d}=}
\newcommand{\I}{\mathbb{I}}

\newcommand{\F}{\mathcal F}
\newcommand{\abs}[1]{\left|#1\right|}

\newcommand{\sizeb}[1]{#1^{*}}
\newcommand{\zb}[1]{#1^{(z)}}
\newcommand{\shb}[1]{{#1}^{\scriptscriptstyle \Box}}

\title{A square bias transformation: properties and applications
\thanks{Research supported by the Russian
Foundation for Basic Research (projects 11-01-00515a, 11-07-00112a,
11-01-12026-ofi-m) and by the grant of the President of Russia
(MK--2256.2012.1).}}
\author{Irina Shevtsova\thanks{Lomonosov Moscow State University,
Faculty of Computational Mathematics and Cybernetics, Leninskie
Gory, GSP-1, Moscow, 119991, Russia; Institute for Informatics
Problems of the Russian Academy of Sciences; e-mail:
ishevtsova@cs.msu.su}}
\date{}

\begin{document}

\maketitle

\begin{abstract}
The properties of the square bias transformation are studied, in
particular, the precise moment-type estimate for the $L_1$-metric
between the transformed and the original distributions is proved, a
relation between their characteristic functions is found. As a
corollary, some new moment-type estimates for the proximity of
arbitrary characteristic function with zero mean and finite third
moment to the normal one with zero mean and the same variance are
proved involving the double integrals of the square- and zero- bias
transformations.
\end{abstract}

\smallskip


\smallskip


\smallskip

\section{Introduction}

Let $X$ be a random variable (r.v.) with the distribution function
(d.f.) $F(x)=\Prob(X<x)$, $x\in\R,$ and the characteristic function
(ch.f.)
$$
f(t)=\E e^{itX}=\int_{-\infty}^\infty e^{itx}\,dF(x),\ t\in\R,
$$
which is the Fourier-Stieltjes transform of the d.f. $F(x)$. As is
well known, if $X$ is nonnegative with $0<\E X<\infty$, then
\begin{equation}\label{SizeBiasChF}
\frac{f'(t)}{f'(0)}, \quad t\in\R,
\end{equation}
is a ch.f., and if $0<\E X^2<\infty$, then
$$
\frac{f'(t)-f'(0)}{tf''(0)},\quad \frac{f''(t)}{f''(0)}
$$
are ch.f.'s as well (see, e.g.,~\cite[Theorem 12.2.5]{Lukacs1970}).
The probability transformation given by~\eqref{SizeBiasChF} is
called the \textit{$X$-size bias transformation}. By a
transformation of a random variable we mean that of its
distribution. The $X$-size bias transformation was introduced by
Goldstein and Rinott~\cite{GoldsteinRinott1996} for the purpose of
estimation of the accuracy of the multivariate normal approximation
to nonnegative random vectors under conditions of local dependence
by Stein's method. Namely, in~\cite{GoldsteinRinott1996} an almost
surely nonnegative r.v. $\sizeb{X}$ with $0<\E X<\infty$ is said to
have the $X$-size biased distribution if
\begin{equation}\label{SizeBiasDF}
d\Prob(\sizeb{X}<x)=\frac{x}{\E X}\,dF(x),\quad x\in\R.
\end{equation}
It is easy to see that the distribution given by~\eqref{SizeBiasDF}
has the ch.f. given by~\eqref{SizeBiasChF}, hence, by virtue of the
uniqueness theorem, definitions~\eqref{SizeBiasChF}
and~\eqref{SizeBiasDF} are equivalent. In the same paper Goldstein
and Rinott also noticed that the distribution of $\sizeb{X}$ may be
characterized by the relation
$$
\E XG(X)=\E X\E G(\sizeb{X}),
$$
which should hold for all functions $G$ such that $\E XG(X)<\infty$.

As regards the second transformation, if $\E X=0$, then the
distribution given by the ch.f.
$$
\frac{f'(t)-f'(0)}{tf''(0)} =-\frac1{\sigma^2}\cdot\frac{f'(t)}{t},
$$
where $\sigma^2=\E X^2>0$, is called the \textit{$X$-zero biased
distribution}. This definition was introduced by Goldstein and
Reinert in~\cite{GoldsteinReinert1997} in an equivalent form for the
purpose of generalization of the size bias transformation to r.v.'s
taking both positive and negative values and was inspired by the
characteristic property of the mean zero normal distribution as the
unique fixed point of the zero bias transformation. Namely,
in~\cite{GoldsteinReinert1997} a r.v. $\zb{X}$ is said to have the
$X$-zero biased distribution if $\E X=0$ and
$$
\E XG(X)=\sigma^2\E G'(\zb{X})
$$
for all absolutely continuous functions $G$ for which $\E XG(X)$
exists. The zero biased transformation possesses the following
elementary properties (most of them are noticed/proved
in~\cite{GoldsteinReinert1997}):

\begin{enumerate}

\item The zero biased distribution is absolutely continuous and
unimodal about zero with the probability density function
$$
p(x)=\sigma^{-2}\E X\I(X>x),\quad x\in\R,
$$
and the ch.f.
$$
\E e^{it\zb{X}}\equiv -\frac1{\sigma^2}\cdot\frac{f'(t)}{t}, \quad
t\in\R.
$$

\item $\zb{X}\eqd X$ if and only if $X$ has the normal distribution
with zero mean~\cite{Stein1981,GoldsteinReinert1997}.

\item The zero biased transformation preserves
symmetry.

\item $\sigma^2\E(\zb{X})^n=\E X^{n+2}/(n+1)$ for $n\in\N$, in
particular, $\sigma^2\E\zb{X}=0.5\E X^3$.

\item If $X=Y_1+\ldots Y_n$, where $Y_1,\ldots,Y_n$ are
independent r.v.'s with zero means and $\E Y_j^2=\sigma_j^2>0$ so
that $\sigma_1^2+\ldots+\sigma_n^2=\sigma^2$, then $\zb{X}=
X_I+\zb{Y_I},$ where $I$ is a random index independent of
$Y_1,\ldots,Y_n$ with the distribution
$\Prob(I=i)=\sigma_i^2/\sigma^2$, $i=1,\ldots,n,$ and
$X_i=X-Y_i=\sum_{j\neq i}Y_j$.

\item
The following non-trivial estimate was proved in 2009 independently
by Goldstein~\cite{Goldstein2009} and
Tyurin~\cite{Tyurin2009DAN,Tyurin2009arxiv}:
\begin{equation}\label{ZeroBiasL1}
L_1(X,\zb{X})\le \frac{\E|X|^3}{2\sigma^2},
\end{equation}
$L_1(X,Y)$ being the $L_1$-distance between the r.v.'s $X$ and $Y$,
$$
L_1(X,Y)= \inf\left\{ \E|X'-Y'|\colon X'\eqd X,\ Y'\eqd
Y\right\},\quad \E|X|<\infty,\ \E|Y|<\infty.
$$
\end{enumerate}

As regards the third transformation, given by the characteristic
function
\begin{equation}\label{ShapeBiasChF}
\shb{f}(t)\equiv\frac{f''(t)}{f''(0)}=-\frac{f''(t)}{\sigma^2},\
t\in\R,
\end{equation}
where $\sigma^2\equiv \E X^2\in(0,\infty)$, it is called the
\textit{$X$-square bias transformation}. It is easy to see that a
r.v. $\shb{X}$ has the ch.f. $\shb{f}(t)$ if and only if
\begin{equation}\label{ShapeBiasCharacterizationG}
\E X^2G(X)=\sigma^2\E G(\shb{X})
\end{equation}
for all functions $G$ such that $\E X^2|G(X)|<\infty$. In 2007
L.\,Goldstein~\cite{Goldstein2007} called the distribution of a r.v.
$\shb{X}$ satisfying~\eqref{ShapeBiasCharacterizationG} the
$X$-square biased distribution. In 2011 L.\,Chen, L.\,Goldstein and
Q.-M.\,Shao~\cite[Proposition~2.3]{ChenGoldsteinShao2011} proved the
following relation between the distribution of the zero biased
$\zb{X}$ and square biased $\shb{X}$ distributions of a symmetric
r.v. $X$:
$$
\zb{X}\eqd U\shb{X},
$$
where the r.v.'s $U$, $\shb{X}$ are independent, $U$ having uniform
distribution on $[-1,1]$. Taking into account that the square bias
transformation preserves symmetry (see below), the latest relation

Later in~\cite{PekozRollinRoss2013} it was noticed that
$\shb{X}=\sizeb{(\sizeb{X})}$, and the distribution of $\shb{X}$
obtained the second name: $X$-double size bias distribution. In the
same paper the following characterization of the modulus of the
normally distributed r.v. was proved: $X\eqd U_1\shb{X}$, where
$U_1$ has uniform distribution on $[0,1]$ and is independent of
$\shb{X}$, if and only if $X\eqd|Z|$, where $Z$ has the standard
normal distribution.

It is easy to see that the square bias transformation possesses the
following elementary properties.

\begin{enumerate}

\item A r.v. $\shb{X}$ has the $X$-square biased distribution if and
only if its d.f. $\shb{F}(x)$ satisfies
\begin{equation}\label{ShapeBiasDF}
d\shb{F}(x)=\frac{x^2}{\sigma^2}\,dF(x),\ x\in\R.
\end{equation}

\item $\shb{X}\eqd X$ if and only if $\Prob(|X|=\sigma)=1$, i.\,e. any
Bernoulli distribution with symmetric atoms is a fixed point of the
square bias transformation. This can be verified by noticing that
the solution of the corresponding linear homogeneous differential
equation $f''(t)+\sigma^2f(t)=0$ of the second order with the
initial condition $f(0)=1$ has the form $f(t)=pe^{i\sigma
t}+(1-p)e^{-i\sigma t}$, $p\in\R$, being a ch.f. if and only if
$p\in[0,1]$.

\item
Square bias transformation preserves symmetry. Indeed, if the r.v.
$X$ has a symmetric distribution, then it's ch.f. $f(t)$ is even,
and hence, $\shb{f}(-t)=f''(-t)/f''(0)=f''(t)/f''(0)=\shb{f}(t)$,
i.e. the distribution of $\shb{X}$ is symmetric as well.

\item
$(\shb{X})^2\eqd\sizeb{(X^2)}$, where $\sizeb{(X^2)}$ has the
$X^2$-size biased distribution.

\item $\shb{(cX)}=c\shb{X}$ for any constant $c\in\R$.

\item
$\sigma^2\E (\shb{X})^n=\E X^{n+2},$ $n\in\N$, $\sigma^2\E
|\shb{X}|^r=\E|X|^{r+2},$ $r>0,$  in particular, $\sigma^2\E
\shb{X}=\E X^3,$ $\sigma^2\E|\shb{X}|=\E|X|^3$.

\end{enumerate}

Moreover, the following estimate for the $L_1$-distance between the
distributions of $X$ and $\shb{X}$ will be proved in this paper.

\begin{theorem}\label{ThL1(X,X-shape_bias)}
If $\E X=0,$ $\E X^2=1,$ and $\E|X|^3<\infty,$ then
$$
L_1(X,\shb{X})\le \E|X|^3,
$$
moreover, for any $\varepsilon>0$ there exists a distribution of a
r.v. $X$ concentrated in two points, such that $\E X=0,$ $\E X^2=1,$
$\E|X|^3<\infty,$ and
$$
L_1(X,\shb{X})> (1-\varepsilon)\E|X|^3.
$$
\end{theorem}

The existence of the square bias transformation follows from the
earlier result of~\cite{Lukacs1970} mentioned above. Moreover, in
2005, Goldstein and Reinert~\cite{GoldsteinReinert2005} proved the
existence of a class of transformations of probability distributions
that are characterized by equations
like~\eqref{ShapeBiasCharacterizationG}. Namely, the authors
described a class of measurable functions $T\colon\R\to\R$ that
provide the existence and uniqueness of the distribution of a random
variable $X^{(T)}$ such that
$$
\E T(X)G(X)=\E G^{(m)}(X^{(T)})\cdot\frac{\E X^mT(X)}{m!}
$$
for all $m$ times differentiable functions $G\colon\R\to\R$ with
$\E|T(X)G(X)|<\infty$. The authors of~\cite{GoldsteinReinert2005}
also noticed that this class includes the zero- and size- bias
transformations respectively with $m=1$, $T(x)=x$ and $m=0$,
$T(x)=x^{+}$. L.\,Goldstein~\cite{Goldstein2007} noticed that this
class also includes the square bias transformation (with $m=0$,
$T(x)=x^2$). However, up till now the properties of the square
biased transformation have not been studied, in particular, the
characteristic function of the square biased distribution and the
estimate for $L_1(X,\shb{X})$ are established in this paper for the
first time.

\section{Motivation and applications}\label{SectionShapeBiasMotivation}

The zero bias transformation gives an opportunity to construct an
integral estimate for the proximity of a ch.f. with zero mean to the
normal one with the same variance in terms of the proximity of the
corresponding zero biased distribution to the original one, which
might be sharper than non-integral estimates based on the Taylor
formula in the neighborhood of zero. Namely, for the sake of
convenience put $\sigma^2=1$ implying $\beta_3\ge1$ by the Lyapounov
inequality. Then using the elementary relations
\begin{multline}\label{ChFDiffInt1}
f(t)-e^{-t^2/2}= e^{-t^2/2}\int_0^t\big(f(u)e^{u^2/2}-1\big)'du=
e^{-t^2/2} \int_0^t\big(f'(u)+\sigma^2uf(u)\big)e^{u^2/2}du=
\\
= e^{-t^2/2}\int_0^t\big(\E e^{iuX}-\E e^{iu\zb{X}}\big)
ue^{u^2/2}du,
\end{multline}
and the estimate for the difference of arbitrary ch.f.'s with finite
first moments due to Korolev and
Shevtsova~\cite{KorolevShevtsova2010SAJ}:
\begin{equation}\label{ChFDiffSinKSh}
\abs{\E e^{itX}-\E e^{itY}}\le 2\sin\Big(
\frac{|t|}{2}\,L_1(X,Y)\wedge\frac\pi2\Big),\quad \E|X|,
\E|Y|<\infty,\quad t\in\R,
\end{equation}
where $a\wedge b\equiv\min\{a,b\}$, $a,b\in\R$, it is not difficult
to conclude that
$$
r(t)\equiv\abs{f(t)-e^{-t^2/2}}\le e^{-t^2/2}\int_0^{|t|}\abs{\E
e^{iuX}-\E e^{iu\zb{X}}} ue^{u^2/2}du\le
$$
$$
\le 2e^{-t^2/2}\int_0^{|t|}\sin\Big(
\frac{u}{2}L_1(X,\zb{X})\wedge\frac\pi2\Big) ue^{u^2/2}du,\quad
t\in\R.
$$
Finally, applying inequality~\eqref{ZeroBiasL1} to estimate
$L_1(X,\zb{X})$ one obtain
\begin{equation}\label{ChFDiffEstimZeroBiasInt1_n=1}
r(t)\le 2e^{-t^2/2}\int_0^{|t|}\sin\Big(
\frac{\beta_3u}{4}\wedge\frac\pi2\Big) ue^{u^2/2}du,\quad t\in\R,
\end{equation}
for any r.v. $X$ with $\E X=0$, $\E X^2=1$,
$\E|X|^3=\beta_3<\infty$.
Estimate~\eqref{ChFDiffEstimZeroBiasInt1_n=1} is exact as $t\to0$,
since, as is well known, $r(t)\sim\beta_3|t|^3/6$,
and~\eqref{ChFDiffEstimZeroBiasInt1_n=1} implies that for all
$\beta_3\ge1$ such that $\beta_3|t|\le\pi/2$ we have
$$
r(t)\le 2\int_0^{|t|}u\sin\Big( \frac{\beta_3u}{4}\Big)du
=\frac{32}{\beta_3^2}\Big(\sin\frac{\beta_3|t|}{4}
-\frac{\beta_3|t|}{4}\cos\frac{\beta_3|t|}{4}\Big)<
\frac{\beta_3|t|^3}{6},
$$
with the least possible factor $1/6$. However,
estimate~\eqref{ChFDiffEstimZeroBiasInt1_n=1} is always sharper than
the power-type estimate $r(t)\le\beta_3|t|^3/6$ especially for
moderate (separated from zero) values of $\beta_3|t|$. Note that
$\beta_3|t|$ can be separated from zero for large enough values of
$\beta_3\ge1$ even if $|t|$ is small. Thus, the estimates for $r(t)$
of an integral \eqref{ChFDiffEstimZeroBiasInt1_n=1}-type form play
an important role in the construction of the least possible upper
moment-type bounds of the accuracy of the normal approximation which
should be uniform in some classes of distributions, especially if in
these classes extremal distributions have large third absolute
moments. This situation is typical, for example, for the problem of
optimization of the absolute constants in the Berry--Esseen-type
inequalities with an improved structure (see
\cite{KorolevShevtsova2010DAN1,KorolevShevtsova2009,KorolevShevtsova2010SAJ,KorolevShevtsova2010TVP,Shevtsova2011arxiv}
where a smoothing inequality is applied with the subsequent
estimation of the difference $|f_n(t)-e^{-t^2/2}|$, $f_n(t)$ being
the ch.f. of the normalized sum of independent random variables, in
terms of the difference $|f(t)-e^{-t^2/2}|$, $f(t)$ being the ch.f.
of a single r.v.) and in its non-uniform analogues for sums of
independent r.v.'s that use the Berry--Esseen inequality with an
improved structure
(see~\cite{Gavrilenko2011,NefedovaShevtsova2012,GrigorievaPopov2012}),
as well as in the moment-type estimates of the rate of convergence
in limit theorems for compound and mixed compound Poisson
distributions (where $\beta_3\to\infty$, see
\cite{KorolevShevtsova2010DAN2,KorolevShevtsova2010SAJ,NefedovaShevtsova2011DAN})
which use the Berry--Esseen inequality with an improved structure as
well.

The above reasoning suggests that for the moderate values of
$\beta_3|t|$, estimates for $r(t)$ in the twice-integrated form
might be even sharper than estimates in the once-integrated form
like~\eqref{ChFDiffEstimZeroBiasInt1_n=1}. Since the ch.f. $f(t)$ is
supposed to be differentiable at least twice, it is possible to
continue~\eqref{ChFDiffInt1} as
\begin{multline*}
f(t)-e^{-t^2/2}= e^{-t^2/2} \int_0^te^{u^2/2}\int_0^u
\big(f'(s)+sf(s)\big)'ds\,du=
\\
= e^{-t^2/2} \int_0^te^{u^2/2}\int_0^u
\big(f''(s)+f(s)+sf'(s)\big)ds\,du=
\\
= e^{-t^2/2}\int_0^te^{u^2/2}\int_0^u \Big(\E e^{isX}-\E
e^{is\shb{X}} + sf'(s)\Big) ds\,du,
\end{multline*}
or as
\begin{multline*}
f(t)-e^{-t^2/2}= e^{-t^2/2}
\int_0^t\int_0^u\big(f(s)e^{s^2/2}-1\big)''ds\,du=
\\
= e^{-t^2/2} \int_0^t\int_0^u\Big( f''(s)+f(s)+ sf'(s)
+s\big(f'(s)+sf(s)\big) \Big)e^{s^2/2}ds\,du=
\\
= e^{-t^2/2}\int_0^t\int_0^u \Big(\E e^{isX}-\E e^{is\shb{X}} +
sf'(s)+s^2\big(\E e^{isX}-\E e^{is\zb{X}}\big)\Big) e^{s^2/2}ds\,du.
\end{multline*}
Note that the second estimate contains the additional term
$s^2\big(\E e^{isX}-\E e^{is\zb{X}}\big)$, but the factor
$e^{s^2/2}$ does not exceed the analogous factor $e^{u^2/2}$ in the
first one. However, for all $0\le s\le t$ this additional term
satisfies
$$
g_3(s)\equiv s^2\big|\E e^{isX}-\E e^{is\zb{X}}\big|\le
2s^2\le2t^2=O(t^2), \quad t\to0,
$$
(actually, an even sharper estimate can be obtained, if
inequalities~\eqref{ChFDiffSinKSh} and~\eqref{ZeroBiasL1} are used).
If $\E X=0$, $\E X^2=1$, then for all $0\le s\le t$ we have
$$
g_2(s)=s|f'(s)|= s|\E Xe^{isX}-\E X|\le s\E|X(e^{isX}-1)|\le s^2\E
X^2=s^2\le t^2 =O(t^2),\quad t\to0,
$$
so that
$$
\sup_{0\le s\le t}(g_2(s)+g_3(s))= O(t^2),\quad t\to0,
$$
while the first term
$$
g_1(s)= \abs{\E e^{isX}-\E e^{is\shb{X}}}=\abs{f(s)+f''(s)}
$$
should be equivalent to $\beta_3s$ as $s\to0+$ in order that the
final integrated estimate have the exact order ${\beta_3|t|^3}/6$ as
$t\to0$. Thus, it is $g_1(s)$ that determines the behavior of the
final integral estimate for small values of $s$, and the problem of
construction of the least possible bound for $g_1(s)$ is very
important. Theorem~\ref{ThL1(X,X-shape_bias)} gives an opportunity
to construct such a bound. Namely, the following corollaries hold.

\begin{cor}
Let $X$ be a r.v. with the ch.f. $f(t)$ and $\E X=0$, $\E X^2=1$,
$\beta_3\equiv\E|X|^3<\infty$. Then for all $t\in\R$
$$
\abs{\E e^{isX}-\E e^{is\shb{X}}}\equiv \abs{f(t)+f''(t)} \le
2\sin\Big( \frac{\beta_3|t|}{2}\wedge\frac\pi2\Big).
$$
\end{cor}

\begin{cor}\label{CorChFDiffShapeBiasInt2_n=1}
Let $X$ be a r.v. with the ch.f. $f(t)$ and $\E X=0$, $\E X^2=1$,
$\beta_3\equiv\E|X|^3<\infty$. Then for all $t\in\R$
$$
\abs{f(t)-e^{-t^2/2}}\le e^{-t^2/2}\int_0^{|t|} \bigg(
2\int_0^u\sin\Big(\frac{\beta_3s}{2}\wedge\frac\pi2\Big)ds +
\frac{u^3}{3}\bigg)e^{u^2/2}\wedge
$$
$$
\wedge\int_0^u \Big(2\sin\Big(\frac{\beta_3s}{2}\wedge\frac\pi2\Big)
+ 2s^2\sin\Big( \frac{\beta_3s}{4}\wedge\frac\pi2\Big) +s^2\Big)
e^{s^2/2}ds\,du.
$$
\end{cor}

Note that, as $t\to0+$, the r.-h. sides of the inequalities
presented in corollary~\ref{CorChFDiffShapeBiasInt2_n=1} are
equivalent to
\begin{multline*}
2\int_0^t \int_0^u \sin\Big(\frac{\beta_3s}{2}\Big)ds\,du=
\frac8{\beta_3^2}\Big(\frac{\beta_3t}2-\sin\frac{\beta_3t}2 \Big) <
\\
< \frac{32}{\beta_3^2}\Big(\sin\frac{\beta_3t}{4}
-\frac{\beta_3t}{4}\cos\frac{\beta_3t}{4}\Big)=
2\int_0^{t}u\sin\Big( \frac{\beta_3u}{4}\Big)du,
\end{multline*}
provided that $0\le\beta_3t\le\pi$. Thus, the estimates including
the square bias transformation which are presented in
corollary~\ref{CorChFDiffShapeBiasInt2_n=1} in the twice-integrated
form are sharper as $t\to0$ that the estimates in the
once-integrated form which include the zero bias transformation
only. So, corollary~\ref{CorChFDiffShapeBiasInt2_n=1} plays an
important role in estimation of the rate of convergence in limit
theorems for sums of independent random variables mentioned above.
However, particular application of
corollary~\ref{CorChFDiffShapeBiasInt2_n=1} is the subject of a
separate investigation and will be published elsewhere.

\section{Proof of theorem~\ref{ThL1(X,X-shape_bias)}}

As is known (see, e.g.~\cite[Theorem~1.3.1]{Zolotarev1986}), the
$L_1$-metric can be represented in terms of the $\zeta_1$-metric as
$$
L_1(X,Y)=\zeta_1(X,Y)\equiv \sup\{|\E g(X)-\E g(Y)|\colon
g\in\F_1\},\quad \E|X|<\infty, \E|Y|<\infty,
$$
where $\F_1$ is the set of all real-valued functions on $\R$ such
that $\sup_{x\neq y}|g(x)-g(y)|/|x-y|\le1$. Since for any function
$g\in\F_1$ we also have $(-g)\in\F_1$, we conclude that the modulus
in the definition of $\zeta_1(X,Y)$ can be omitted:
$$
L_1(X,\shb{X})=\zeta_1(X,\shb{X})=\sup\{\E g(X)-\E g(\shb{X})\colon
g\in\F_1\}.
$$

Let $X$ be a r.v. with the d.f. $F(x)$ and $\E X=0$, $\E X^2=1$,
${\E|X|^3<\infty}$, $\shb{X}$ have $X$-square biased distribution,
i.e. the d.f. $\shb{F}(x)$ of the r.v. $\shb{X}$ satisfying the
relation $d\shb{F}(x)=x^2dF(x)$, $x\in\R$. Then
$$
L_1(X,\shb{X})=\sup_{g\in\F_1}(\E g(X)-\E g(\shb{X}))=
\sup_{g\in\F_1}\int_{-\infty}^\infty(1-x^2)g(x)dF(x).
$$
For $g\in\F_1$ denote
$$
J(F,g)=\E((1-X^2)g(X)-|X|^3)=
\int_{-\infty}^\infty((1-x^2)g(x)-|x|^3)dF(x).
$$
Then
$$
L_1(X,\shb{X})-\E|X|^3=\sup_{g\in\F_1}J(F,g),
$$
and the statement of the theorem is equivalent to
$$
\sup_{g\in\F_1}\sup_{F} J(F,g)=0,
$$
where the supremum $\sup_F$ is taken over all d.f.'s $F$ of the r.v.
$X$ satisfying two moment-type conditions: $\E X=0$, $\E X^2=1$. As
it follows from the results
of~\cite{Hoeffding1955,MulhollandRogers1958}, the supremum of a
linear (with respect to the d.f. $F(x)$) functional $J(X,g)$ under
two linear equality-type conditions $\E X=0$, $\E X^2=1$ is attained
at the distributions concentrated in at most three points.

Before passing to checking three- and two-point distributions,
recall that for $L_1$-metric the following representation in terms
of the mean metric holds as well (see,
e.g.~\cite[\S~1.3]{Zolotarev1986}):
$$
L_1(X,Y)=\varkappa(X,Y)\equiv
\int_{-\infty}^\infty|\Prob(X<u)-\Prob(Y<u)|\,du,\quad
\E|X|\vee\E|Y|<\infty,
$$
and hence,
$$
L_1(X,\shb{X}) =  \int_{-\infty}^\infty|F(u)-\shb{F}(u)|\,du =
\int_{-\infty}^\infty|F(u)-\E X^2\I(X<u)|\,du.
$$
$$
L_1(X,\shb{X}) = \int_{-\infty}^\infty|F(u)-\E X^2\I(X<u)|\,du.
$$

Let the r.v. $X$ take two values and satisfy the conditions $\E
X=0$, $\E X^2=1$. Then its distribution should necessarily have the
form
$$
\Prob\left(X=\sqrt{q/p}\right)=p=1-\Prob\left(X=-\sqrt{p/q}\right),
\quad q=1-p\in(0,1).
$$
It is easy to see that $\E X^3=(q-p)/\sqrt{pq}$,
$\E|X|^3=(p^2+q^2)/\sqrt{pq}$. Then
$$
\E X^2\I(X<u)= \left\{
  \begin{array}{ll}
    0, & u\le -\sqrt{p/q}, \\
    p, & -\sqrt{p/q}<u\le\sqrt{q/p}, \\
    1, & \sqrt{q/p}<u,
  \end{array}
\right.
$$
$$
\E X^2\I(X<u)-F(u)= (p-q)\I\big(-\sqrt{p/q}<u\le\sqrt{q/p}\,\big),
$$
and hence
$$
L_1(X,\shb{X})= \int_{-\infty}^\infty
|p-q|\I\big(-\sqrt{p/q}<u\le\sqrt{q/p}\,\big)\,du =
\frac{|p-q|}{\sqrt{pq}}=|\E X^3|\le \E|X|^3
$$
by virtue of the Jensen inequality, thus, the statement of the
theorem holds. Moreover, for any $\eps>0$
$$
\frac{L_1(X,\shb{X})}{\E|X|^3}
=\frac{|1-2p|}{1-2p+2p^2}\ge1-2p>1-\eps
$$
for all $0<p<\eps/2$.

Now consider a r.v. $X$ taking exactly three values. Note that
$$
\sup\Big\{\frac{L_1(X,\shb{X})}{\E|X|^3}\colon \E X=0,\ \E
X^2=1\Big\}= \sup\Big\{\frac{L_1(X,\shb{X})\sigma^2}{\E|X|^3}\colon
\sigma>0,\ \E X=0,\ \E X^2=\sigma^2\Big\},
$$
where the supremums are taken over three-point distributions of the
r.v. $X$. Let $X$ take values $x,y,z$ with probabilities $p,q,r>0$
respectively, $p+q+r=1$. Without loss of generality one can assume
that $x<y\le0<z$. From the conditions $\E X=0,\ \E X^2=\sigma^2$ we
find that
$$
p=\frac{\sigma^2+yz}{(z-x)(y-x)},\quad
q=-\frac{\sigma^2+xz}{(z-y)(y-x)},\quad
r=\frac{\sigma^2+xy}{(z-x)(z-y)},\quad -yz<\sigma^2<-xz.
$$
For all $u\in\R$ we have
$$
\E X^2\I(X<u)= \left\{
  \begin{array}{ll}
    0, & u\le x, \\
    px^2, & x<u\le y, \\
    px^2+qy^2, & y<u\le z, \\
    \sigma^2, & z<u,
  \end{array}
\right.
$$
$$
\sigma^{-2}\E X^2\I(X<u)-F(u)= \left\{
  \begin{array}{ll}
    0, & u\le x, \\
    p(x^2/\sigma^2-1), & x<u\le y, \\
    (px^2+qy^2)/\sigma^2-p-q, & y<u\le z, \\
    0, & z<u.
  \end{array}
\right.
$$
Noticing that $(px^2+qy^2)/\sigma^2-p-q=
(\sigma^2-rz^2)/\sigma^2-1+r=r(1-z^2/\sigma^2)$, we obtain
$$
L_1(X,\shb{X}) = \int_{-\infty}^\infty|F(u)-\sigma^{-2}\E
X^2\I(X<u)|\,du= p\Big|\frac{x^2}{\sigma^2}-1\Big|(y-x)+
r\Big|1-\frac{z^2}{\sigma^2}\Big|(z-y).
$$
Consider the function
\begin{multline*}
{L_1(X,\shb{X})\sigma^2}-{\E|X|^3}=
p(y-x)|x^2-\sigma^2|+r(z-y)|z^2-\sigma^2|+px^3+qy^3-rz^3=
\\
=\frac1{z-x}\Big(\abs{x^2-\sigma^2}(\sigma^2+yz) +
\abs{z^2-\sigma^2}(\sigma^2+xy)-\frac{2z^3(\sigma^2+xz)}{z-y}+
\\
+ \sigma^2(z^2-x^2-xy+yz) +
 xyz(z-x)+2xz^3\Big) \equiv g(x,y,z,\sigma^2).
\end{multline*}
The statement of the theorem is equivalent to $\sup
g(x,y,z,\sigma^2)=0$, where the sumpremum is taken over all
$\sigma^2>0$, $x<y\le0<z$ such that $-yz<\sigma^2<-xz$. Note that it
suffices to consider only $\sigma^2<\max\{x^2,z^2\}$, since the
opposite inequality (with $q>0$) implies that $\E X^2<\sigma^2$. So,
there are only three possibilities: 1)~$0<\sigma^2<\min\{x^2,z^2\}$,
2)~$x^2\le\sigma^2<z^2$, 3)~$z^2\le\sigma^2<x^2$. Opening the
modules, we notice that $g(x,y,z,\sigma^2)$ is a parabola with
respect to $\sigma^2$ on each of the intervals specified above.
Consider the behavior of $g(x,y,z,\sigma^2)$ on each of these
intervals.

\begin{enumerate}

\item
$0<\sigma^2<\min\{x^2,z^2\}$, then necessarily $\sigma^2<-xz$ and
$$
g(x,y,z,\sigma^2)  =
-\frac{2(\sigma^2+xy)(yz^2+\sigma^2(z-y))}{(z-y)(z-x)}.
$$
The coefficient $-2/(z-x)$ at $\sigma^4$ is negative, thus the
branches of this parabola with respect to $\sigma^2$ look down and
the maximum value of the function $g(z,y,z,\sigma^2)$ within the
interval $0<\sigma^2<\min\{x^2,z^2\}$ is attained either at the
vertex
$$
\sigma^2_*=-\frac{y(z^2+xz-xy)}{2(z-y)},
$$
if $\sigma^2_*>-yz$, or at the point $\sigma^2\to-yz+0$, if
$\sigma^2_*\le-yz$. We have
$$
\sigma^2_*+yz= \frac{y(z(z-x)-y(2z-x))}{2(z-y)}\le 0,
$$
since $x<y\le0<z$, with the equality attained if and only if $y=0$.
Thus, the supremum is attained as $\sigma^2\to-yz+0$, which implies
that $p\to0$ and reduces the problem to checking two-point
distributions considered above.

\item
$z^2\le\sigma^2< x^2$, then
$$
g(x,y,z,\sigma^2)= -\frac{2z^3(\sigma^2+xy)}{(z-x)(z-y)}=-2rz^3<0
$$
by virtue of the conditions $r>0$, $z>0$.

\item
$x^2\le\sigma^2< z^2$, then the function
$$
g(x,y,z,\sigma^2)= 2\frac{-\sigma^2(x^2(z-y)+y^2(z-x)+xyz)
+xyz(xy-xz-yz)}{(z-x)(z-y)}
$$
is linear and decreases monotonically in $\sigma^2$, since
$x^2(z-y)+y^2(z-x)+xyz>0$. Thus, if $x^2\le-yz$, then the supremum
of $g(x,y,z,\sigma^2)$ is supplied by $\sigma^2\to-yz+0$, which
reduces the problem to checking two-point distributions considered
above. If $x^2>-yz$, then the supremum of $g(x,y,z,\sigma^2)$ is
attained at $\sigma^2=x^2$. With this value of $\sigma^2$ we have
$$
g(x,y,z,x^2)=-\frac{2x(x+y)(x^2z-yx^2+yz^2)}{(z-x)(z-y)}.
$$
$$
\frac{\partial}{\partial y}\,
g(x,y,z,x^2)=-\frac{2x(x+z)(yz(2z-y)+x(z-y)^2)}{(z-x)(z-y)^2}<0,
$$
since $x<-z$ and $x<y\le0<z$. Thus, the supremum of $g(x,y,z,x^2)$
over all $y$ such that $x<y\le0$ and $-yz<x^2=\sigma^2$ is supplied
by $y\to\max\{x,-x^2/z\}+0=-x^2/z+0$, i.\,e., $p\to0$, which reduces
the problem to checking two-point distributions considered above.
Thus, the theorem is completely proved.
\end{enumerate}


\end{document}